\documentstyle[12pt]{article}
\newcommand{\const}{\mathop{\rm const}\limits}

\newcommand{\mod}{\mathop{\rm mod}\limits}

\newcommand{\Law}{\mathop{\rm Law}\limits}

\newcommand{\Comp}{\mathop{\rm Comp}\limits}

\newcommand{\Dom}{\mathop{\rm Dom}\limits}

\newcommand{\Ent}{\mathop{\rm Ent}\limits}

\begin{document}

 \begin{center}

{\bf Each Random Variable in separable Banach Space }\\

\vspace{3mm}

{\bf belongs to the Domain of Definition of some }\\

\vspace{3mm}

{\bf  inverse to compact linear non-random operator.  }\\

\vspace{5mm}

 $ {\bf E.Ostrovsky^a, \ \ L.Sirota^b } $ \\

\vspace{4mm}

$ ^a $ Corresponding Author. Department of Mathematics and computer science.\\
 Bar-Ilan University, 84105, Ramat Gan, Israel.\\

\vspace{4mm}

E-mail:  eugostrovsky@list.ru\\

\vspace{4mm}

$ ^b $  Department of Mathematics and computer science. \\
Bar-Ilan University, 84105, Ramat Gan, Israel.\\

\vspace{4mm}

E-mail:  sirota3@bezeqint.net\\

 \vspace{4mm}

 {\bf Abstract.}\\

\end{center}

  \vspace{3mm}

  We give in this short report a very simple proof that arbitrary random variable with Borelian distribution in
separable Banach space belongs with probability one to a pre-image of some linear compact non-random operator. \\

\vspace{4mm}

{\it  Key words and phrases: }  Law of distribution, Banach and  Linear Topological Spaces, random variables (r.v.),
non-random sequence, support of measure, compact imbedded subspace and compact linear operator, Franklin's functional
orthonormal system and constant. \par

 \vspace{6mm}

\section{ Introduction. Notations. Definitions. Statement of problem. Short history. }

 \vspace{3mm}

 Let \par

{\bf A}. \ $ X $ be arbitrary Separable Banach Space (SBS) with a norm $ || \cdot ||X, $ briefly: $ X = (X, || \cdot ||X), $
 equipped with Borelian sigma-algebra $ B; $ \par

\vspace{4mm}

{\bf B}. \ $ (\Omega, F, {\bf P}) $ be probability triple with expectation {\bf E;} \\

\vspace{3mm}

{\bf C.} $ \xi: \Omega \to X $ be fixed random variable (r.v.) defined on the our triple with the values in
the space $ X :$

$$
{\bf P} (\xi \in X) = 1 \eqno(1.1)
$$
with Borelian distribution $  \mu_{\xi}(\cdot) :$

$$
\mu_{\xi}(A) = {\bf P} (\xi \in A), \ A \in B. \eqno(1.2)
$$

\vspace{3mm}

 We will denote the set of all {\it compact  } linear  continuous operators acting from the Banach space $  X_1  $  with values
in the Banach space $  X_2 $ as $ \Comp(X_1, X_2); $ we agree for  brevity   $ \Comp(X, X):= \Comp(X, X). $  \par
 The domain of definition of the inverse operator  $ U \in \Comp(X) $ will be denoted as ordinary $ \Dom(U^{-1}): $

$$
\Dom(U^{-1}) = \{x, \ x \in X, \ \exists y \in X \ \Rightarrow x = Uy \}.
$$

 It is known, see \cite{Buldygin1},  \cite{Ostrovsky2}, \cite{Ostrovsky3},  that there exists a subspace $  Y $
which dependent only on the distribution  $ \xi, $
compact imbedded in the space $  X, $  i.e. such that the unit ball of the space $  Y  $ is pre-compact set of the space $  X, $
which may  be as a capacity of the support $ \mu_{\xi}: $

$$
 \mu_{\xi}(Y) = {\bf P} (\xi \in Y) = 1.\eqno(1.3)
$$

\vspace{4mm}

 {\bf  Our aim in this report is to generalize the equality (1.3) and consider some applications. } \par

 \vspace{3mm}

\section{Main result.}

 \vspace{3mm}

{\bf Theorem 2.1. } There exists a compact operator $ U \in \Comp(X) $ which dependent only on the distribution $ \mu_{\xi} $
such that the inverse operator $  U^{-1} $ is unique defined on the set $ \Dom(U^{-1}) $ and such that $ \xi \in \Dom(U^{-1}) $
almost everywhere:

$$
{\bf P}(U^{-1}(\xi) \in X) = 1. \eqno(2.1)
$$

\vspace{3mm}

{\bf Proof.}\\

\vspace{3mm}

{\bf A.}  Let us start for beginning from the case when the space $  X  $ coincides with the space $  c_0 $ of the sequences
goes to zero:
$$
 c_0 = \{x\} = \{  \vec{x} \} = \{  x_1, x_2, \ldots, x_n, \ldots \}, \ \lim_{n \to \infty} x_n = 0, \ || \ x \ || \stackrel{def}{=} \sup_n |x_n|.
\eqno(2.2)
$$

 Let the random vector $ \xi = \vec{\xi}  $ belongs to the space $ c_0,  $ so that $ \lim_{n \to \infty} \xi_n = 0 $ a.e.  There
exists a strictly positive non-random sequence $ \epsilon = \vec{\epsilon} = \{\epsilon_1, \epsilon_2, \ldots \}  $
tending to zero:  $ \vec{\epsilon} \in c_0 $ such that

$$
\xi_n = \epsilon_n \cdot \eta_n, \hspace{6mm} \lim_{n \to \infty} \eta_n = 0 \ (\mod {\bf P}), \eqno(2.3)
$$
see  \cite{Kantorovicz1}, chapter 7, section 4.\par

 Define the following linear operator $ U: c_0 \to c_0 $ as follows:

$$
U( \{ x_1, x_2, \ldots, x_n, \ldots   \} ) = \{\epsilon_1 x_1, \epsilon_2 x_2, \ldots, \epsilon_n x_n, \ldots   \}; \eqno(2.4)
$$
then $ U \in \Comp(c_0) $ and $ \vec{\xi} \in \Dom(U^{-1}) \ (\mod {\bf P} ). $\par

\vspace{3mm}

{\bf B.} At the same proof may be used in the case when the space $ X  $ is equal to the space $ c_0(Y),  $ which consists
by definition on all the sequences of the elements of the (non necessary to be separable) Banach space $  Y  $ tending to zero:

$$
c_0(Y) = \{ y_i, \ i = 1,2,\ldots, \ y_i \in Y, \ \lim_{n \to \infty} y_n = 0  \}. \eqno(2.5)
$$

\vspace{3mm}

{\bf C.} The general case may be reduced to the considered in (2.5).\par
 Note first of all that by virtue of the universality  of the space $ C(0,1) $ in the class of all separable Banach spaces
 it is sufficient to consider only the case when the space $  B  $ coincides with the classical space of all continuous functions
defined  on the closed interval $  [0,1]: \   B = C[0,1], $ equipped with  uniform norm $ ||f|| := \max_{t\in [0,1]} |f(t)|, $
see  \cite{Ostrovsky2}.\par
 Further, we will use the Franklin orthonormal sequence of functions $ \phi_i = \phi_i(t), \ t \in [0,1], i = 1,2,\ldots,$ see
 \cite{Franklin1}; the detail investigation of these functions are obtained in \cite{Kaczmarz1}, \cite{Ciesielski1}. Note only
 that these orthonormal functions forms the unconditional basis in the space $  C(0,1). $\par

 Let  $ \xi(t), \ t \in [0,1] $ be continuous  with probability one random process: \\
 $ {\bf P} (\xi(\cdot) \in C[0,1]) = 1.  $  It may be expressed in the uniform convergent  Fourier - Franklin series

$$
\xi(t) = \sum_{i=1}^{\infty} \xi_i \ \phi_i(t), \hspace{7mm} \xi_i = \int_0^1 \xi(t) \ \phi_i(t) \ dt. \eqno(2.6)
$$
so that with probability one

$$
\lim_{n \to \infty} || \xi(t) - \sum_{i=1}^n \xi_i \ \phi_i(t)|| = 0. \eqno(2.7)
$$

 Let $  \{  N(k) \}, \ k = 1,2,\ldots  $ be strictly  increasing integer non-random sequence, $  N(1) = 1. $  We can and will
suppose without loss of generality that there exists a continuous strictly  increasing function $  N_1(x), $ defined on the semi - axis
$ [1,\infty), $  such that $  N_1(x)/(x=k) = N(k). $ We will denote such a function as before for brevity $  N(x). $

 Define a sequence of linear operators   $  Q_k[\cdot](t)  $

$$
Q_k[\xi](t) = \sum_{n=N(k) + 1}^{N(k+1)} \xi_n \ \phi_n(t), \hspace{7mm} \zeta(k) = ||Q_k[\xi]||, \eqno(2.8)
$$
and correspondingly for arbitrary (continuous) numerical function $ f = f(t) $

$$
Q_k[f](t) = \sum_{n=N(k) + 1}^{N(k+1)} f_n \ \phi_n(t), \hspace{7mm} f_i = \int_0^1 f(t) \ \phi_i(t) \ dt.
 \eqno(2.8a)
$$

 The sequence $ \zeta = \{  \zeta(k) \} $ belongs the the space $  c_0 $ almost surely.\par
 In accordance with the item {\bf A} of this section there exists a non-random strictly increasing unbounded sequence
$  w(k) $ such that by the properly of the selected sequence $  \{ N(k) \} $

$$
\sum_k w(k) \zeta(k) < \infty \ (\mod {\bf P}).\eqno(2.9)
$$
 The required operator $  U:  C[0,1] \to C[0,1]  $ may be  constructed  as follows

$$
U^{-1}[f](t) \stackrel{def}{=} \sum_{k = 1}^{\infty} w(k) \  Q_k[f](t). \eqno(2.10)
$$

Indeed,

$$
U^{-1}[\xi](t) := \sum_{k = 1}^{\infty} w(k) \  Q_k[\xi](t) = \sum_{k=1}^{\infty} w(k) \sum_{n=N(k) + 1}^{N(k+1)} \xi_n \ \phi_n(t)=
$$

$$
\sum_{n=1}^{\infty} v(n) \ \xi_n \ \phi_n(t),
$$
where the coefficients $  v(n) $ have a form

$$
v(n) = \sum_{k=1}^{\infty} w(k) \cdot I( N(k) + 1 \le n \le N(k+1) ), \eqno(2.11)
$$
where $ I(\cdot) $ denotes the indicator function. \par
 Denote

$$
\Delta(n) = N^{-1}(n-1) - (N^{-1}(n) - 1). \eqno(2.12)
$$

 Note that  the values $  v(n) $ allows the following bilateral estimates:

$$
v_-(n) :=  w \left(\Ent[N^{-1}(n) - 1 ] \right) \cdot (\Delta(n) - 1) \le v(n) \le
$$

$$
w \left(\Ent[N^{-1}(n - 1) ] + 1 \right) \cdot (\Delta(n) + 1) =:v_+(n), \eqno(2.13)
$$
where $  \Ent(z) $ denotes the integer part of the value $  z. $ \par
 The operator $  U[f]  $ has a form

$$
U[f](t)= \sum_{n=1}^{\infty} \frac{ f_n \ \phi_n(t)}{v(n)} = \int_0^1 R(t,s) \ f(s) ds, \eqno(2.14)
$$
i.e. is linear integral compact ("diagonal")  operator with continuous by means of the choosing of the sequence $  \{N(k)\} $
kernel

$$
R(t,s) = \sum_{n=1}^{\infty} \frac{\phi_n(t) \ \phi_n(s)}{v(n)}.\eqno(2.15)
$$

\vspace{4mm}

{\bf Remark 2.1.} Note that the kernel, i.e. the function of two variables $ R(t,s) $ in (2.15) is in addition
symmetrical and positive definite; therefore the introduced operator $  U  $ is self-conjugated relative the
ordinary inner product

$$
(g_1, g_2) := \int_0^1 g_1(t) \ g_2(t) \ dt.
$$

 \vspace{3mm}

\section{ Retaining of moments.}

 \vspace{3mm}

  Recall that the function $ \Phi: R \to R_+ $ is called Young function, iff it is even, convex, continuous, non - negative,
 $  \Phi(u) = 0 \Leftrightarrow u = 0, $  is strictly increasing such that $ \lim _{u \to \infty} \Phi(u) = \infty. $ \par

  By definition, $ \Phi(\cdot) $ satisfies the so-called $ \Delta_2 $ condition,
 write: $ \Phi \in \Delta_2, $  iff

 $$
 \forall \lambda > 0 \Rightarrow \sup_{u > 0} \left[ \frac{\Phi(\lambda u)}{\Phi(u)}  \right] < \infty. \eqno(3.1)
 $$

 For instance, the functions  $ \Phi_p(u) = |u|^p, \ p = \const > 0  $ and
 $ \Psi_{p,r}(u) = \\ |u|^p ( 1 + \ln(e + |u|)   )^r, \ p = \const > 0, \ r = \const $ belong to the set $  \Delta_2. $ \par

The Young function $ \Psi(\cdot) $ 	 is called weaker than $ \Phi(\cdot),  $
 if for all positive constant $ \lambda; \lambda = \const > 0 $

$$
\lim_{u \to \infty} \left[ \frac{\Psi( \lambda u)}{\Phi(u)}  \right] = 0. \eqno(3.2)
$$

Notation: $ \Psi << \Phi. $\par

\vspace{3mm}

 Using for us facts about Orlicz's spaces are presented in monographs \cite{Krasnoselsky1},
\cite{Rao1}, \cite{Rao2}; in this report all the Orlicz spaces are constructed over our probability space. \par

\vspace{4mm}

{\bf Theorem 3.1.} Let $ \xi $ be a random vector with values in the separable Banach space $ (X, || \ \cdot \ || = \ || \ \cdot \ || X) $
such that for some Young function  $ \Phi(\cdot) $ belonging to the class $  \Delta_2 $

$$
{\bf E} \Phi(|| \ \xi \ ||X ) < \infty, \eqno(3.3)
$$
or equally such that the r.v. $ || \ \xi \ ||X  $ belongs to the Orlicz space $  L(\Phi). $ \par

 There exists a compact linear operator $ V: X \to X $ depending only on the distribution $ \Law(\xi) $ and on the function $ \Phi(\cdot) $
such that $ \xi \in \Dom(V^{-1}) $ and moreover

$$
{\bf E} \Phi(|| \ V^{-1}[\xi] \ ||X) < \infty. \eqno(3.4)
$$
or equally such that the r.v. $ || \ V^{-1}\xi \ ||X  $ belongs to the Orlicz space $  L(\Phi). $ \par

\vspace{3mm}

{\bf Proof.} It is sufficient to consider  as before  the case when $  X = C[0,1]. $  On the other words,
let as in the second section  $ \xi(t), \ t \in [0,1] $ be continuous  with probability one random process:
 $ {\bf P} (\xi(\cdot) \in C[0,1] = 1),  $ and suppose in addition

$$
{\bf E} \Phi( \sup_{t \in [0,1]} |\xi(t)|) = {\bf E} \Phi( ||\xi(\cdot)||) < \infty. \eqno(3.5)
$$

  It may be expressed in the uniform convergent  Fourier-Franklin series

$$
\xi(t) = \sum_{i=1}^{\infty} \xi_i \ \phi_i(t), \hspace{7mm} \xi_i = \int_0^1 \xi(t) \ \phi_i(t) \ dt, \eqno(3.6)
$$
and we denote

$$
F_n[f](t) = \sum_{i=1}^n f_i \ \phi_i(t),  \hspace{7mm} f(\cdot) \in C[0,1]. \eqno(3.7)
$$

  It follows from the uniform boundedness principle that there is an absolute constant $  C_F $  (Franklin's constant)
such that

$$
\sup_n ||F_n[f]|| \le C_F \ ||f||, \eqno(3.8)
$$
so that with probability one

$$
\lim_{n \to \infty} || \xi(\cdot) - F_n[\xi](\cdot)|| = 0.
$$

 Further, as long as $ || \xi(\cdot) - F_n[\xi](\cdot)|| \le (C_F + 1) \ ||\xi||,   $ and the function $ \Phi(\cdot)  $
is Young function, it follows from  dominated convergence theorem

$$
\lim_{n \to \infty}{\bf E} \Phi(  ||\xi - F_n[\xi]||_c  ) = 0. \eqno(3.9)
$$

 In the language of the theory of Orlicz spaces (see \cite{Krasnoselsky1}, \cite{Rao1}, \cite{Rao2} ) the equality (3.9)
denotes the {\it moment convergence,} or convergence {\it in mean} the variable \\ $ {\bf E} \Phi(||\xi - F_n[\xi]||_c)  $
to zero as $ n \to \infty.  $\par
 But the function $  \Phi $ satisfies the $ \Delta_2 $ condition, the equality (3.9) means also the following Orlicz
norm convergence:

$$
\lim_{n \to \infty} || \hspace{4mm} || \ \xi - F_n[\xi] \ ||_c \hspace{4mm} ||L(\Phi)   = 0, \eqno(3.10)
$$
where for the sake of definiteness $ ||f||_c $ denotes  the uniform norm $ ||f|| $    and $ ||  \eta ||L(\Phi) $ denotes the Orlicz's norm
relative the Young function $ \Phi $ of the r.v. $ \eta  $ defined on our probability space. \par

\vspace{3mm}

 Let $  \{  M(k) \}, \ k = 1,2,\ldots  $ be strictly  increasing integer non - random sequence, $  M(1) = 1. $  We can and will
suppose without loss of generality that there exists a continuous strictly  increasing function $  M_1(x), $ defined on the semi - axis
$ [1,\infty), $  such that $  N_1(x)/(x=k) = M(k). $ We will denote such a function as before for brevity $  M(x). $

 Define a following sequence of linear operators   $  R_k[\cdot](t)  $

$$
R_k[\xi](t) = \sum_{n=M(k) + 1}^{M(k+1)} \xi_n \ \phi_n(t), \hspace{7mm} \nu(k) = ||R_k[\xi]||_c, \eqno(3.11)
$$
and correspondingly for arbitrary (continuous) numerical function $ f = f(t) $

$$
R_k[f](t) = \sum_{n=M(k) + 1}^{M(k+1)} f_n \ \phi_n(t), \hspace{7mm} f_i = \int_0^1 f(t) \ \phi_i(t) \ dt.
 \eqno(3.12)
$$

 The sequence $ \{  M(k) \ \} $ may be  picked such that

$$
|| \ \nu(k) \ ||L(\Phi) \le 4^{-k} || \ ||\xi(\cdot)||_c \ || L(\Phi). \eqno(3.12)
$$

 The required operator $ V: C[0,1] \to C[0,1] $ may be constructed as follows

$$
V^{-1}[f](t) := \sum_k 2^k R_k[f](t). \eqno(3.13)
$$
 Indeed,

$$
|| \  || \ V^{-1}[\xi]  \ ||_c \ || L(\Phi) \le \sum_k 2^{-k} \ || \  ||\xi(\cdot)||_c \ || L(\Phi) \le || \ ||\xi(\cdot)||_c \ || L(\Phi), \eqno(3.14)
$$
Q.E.D. \par

\vspace{4mm}

{\bf Theorem 3.2.} Let $ \xi $ be again a random vector with values in the separable Banach space $ (X, || \ \cdot \ || = \ || \ \cdot \ || X) $
such that for some Young function  $ \Phi(\cdot) $ not necessary  to be from the set $  \Delta_2 $

$$
{\bf E} \Phi(|| \ \xi \ ||X ) < \infty, \eqno(3.15)
$$
or equally such that the r.v. $ || \ \xi \ ||X  $ belongs to the Orlicz space $  L(\Phi). $ \par
Let also 	$  \Psi(\cdot) $ be other arbitrary Young  function weaker than $  \Phi: \ \Psi <<  \Phi. $
 There exists a compact linear operator $ V: X \to X $ depending only  on the distribution $ \Law(\xi) $ and on both the functions $ \Psi, \ \Phi(\cdot) $
such that $ \xi \in \Dom(V^{-1}) $ and moreover

$$
{\bf E} \Psi(|| \ V^{-1}[\xi] \ ||X) < \infty, \eqno(3.16)
$$
or equally such that the r.v. $ || \ V^{-1}\xi \ ||X  $ belongs to the Orlicz space $  L(\Psi). $ \par

\vspace{3mm}

{\bf Proof} is at the same as before,   except the passing from the equality (3.9) to (3.10):

$$
\lim_{n \to \infty}{\bf E} \Phi(  ||\xi - F_n[\xi]||_c  ) = 0
$$
denotes

$$
\lim_{n \to \infty}  || \ ||\xi - F_n[\xi]||_c  \ ||L(\Psi) = 0.
$$

 \vspace{3mm}

\section{ Applications to the CLT in Banach spaces.}

 \vspace{3mm}

 Let again $ \eta  $ be a (Borelian) r.v. with values in the separable Banach space $  X  $ equipped with norm  $ || \ \cdot \ ||X   = || \cdot \ ||. $
Let $ \eta_i, \ i = 1,2,\ldots $ be independent copies of $  \eta, $  may be defined  on some another sufficiently rich  probability space. Denote

$$
S(n) = n^{-1/2} \sum_{i=1}^n \eta_i. \eqno(4.1)
$$

 By definition, the r.v. $ \eta $ is said to satisfy Central Limit Theorem (CLT) in the space $  X,  $  write $ \eta \in CLT, $
 if the sequence $  S(n) $ converges
weakly (in distribution) as $  n \to \infty $  to some Gaussian distributed in this space r.v. $  S(\infty). $\par

 Evidently, $ \eta, \ S(\infty) $ are mean zero, have finite weak second moment,  have identical covariation operator. Moreover,
r.v. $ \eta $  is {\it pregaussian}. This imply by definition that $ {\bf P} ( S(\infty) \in X) = 1. $ \par

 There are   many works about CLT in different Banach spaces, see, e.g.
\cite{Dudley1} \ - \ \cite{Ratchkauskas3} and reference therein. The majority of aforementioned  works are devoted to the finding of necessary
or sufficient conditions for these Theorem. \par

 M.Ledoux and M.Talagrand in the book \cite{Ledoux1}, p.294-296, Lemma 10.1 proved that if the r.v. $ \eta $ satisfies CLT in the space $  X, $
then

$$
\lim_{t \to \infty} t^2 {\bf P} (||\eta|| > t) = 0.
$$
 As a consequence: for every values $ p \in (0,2) \ \Rightarrow {\bf E} ||\eta||^p < \infty. $\par
 Moreover, if the r.v. $ \eta $ satisfies CLT in the space $  X, $ then

$$
\lim_{t \to \infty} t^2 \sup_n {\bf P} (||S(n)|| > t) = 0,
$$

$$
p \in (0,2) \ \Rightarrow  \sup_n{\bf E} ||S(n)||^p < \infty.
$$

M.Ledoux and M.Talagrand \cite{Ledoux1}, p.296-302 introduced and investigated also the following norm for the r.v. $  \eta: $

$$
||\eta||CLT(X) \stackrel{def}{=} \sup_n ||S(n)||X, \eqno(4.2)
$$
and proved in particular that $  \eta \in CLT(X) \ \Rightarrow ||\eta||CLT < \infty,  $ and that the functional
$ \eta \to  ||\eta||CLT   $ is actually  (complete) norm on the linear space of all r.v. satisfying CLT in the space $  X. $ \par
 On the other words, the finiteness of the norm $ ||\eta||CLT(X)  $ is {\it necessary} condition for the conclusion $ \eta \in CLT(X), $
but it is known that this condition is'nt {\it sufficient.}  \par

\vspace{3mm}

{\bf Proposition 4.1.}  Let the centered r.v. $  \eta $ with values in the space $ X $ is pregaussian. Suppose in addition there exists
a compact linear operator $  U: X \to X $ such that

$$
\sup_n || U^{-1} S(n)  || = K(\Law (\eta)) < \infty. \eqno(4.3)
$$
 Then $ \eta \in CLT(X). $ \par

 \vspace{3mm}

 {\bf Proof. } The convergence of finite-dimensional (cylindrical) distributions of $ S(n) $ to ones for $ S(\infty) $ is obvious.
Further, we have  using  Tchebychev's-Markov  inequality introducing the pre-compact set

$$
Z(t) = U\{B_X(t)\}, \ B_X(t) := \{ x, \ x \in X, \ ||x||X \le t \}:
$$

$$
{\bf P} (S(n) \notin  Z(t)) \le \frac{ K(\Law (\eta))}{t}\le \epsilon \eqno(4.4)
$$
for sufficiently greatest values $  t. $ Thus, the sequence of distributions of the variables $ \Law(S(n)) $ is tight in this space. \par
This completes the proof of proposition 4.1. \par

 The examples of applying of this approach, i.e. by means of assertion (4.1),  may be found in the articles
\cite{Ostrovsky100}-\cite{Ostrovsky103}. \par

 \vspace{3mm}

\section{ Concluding remarks.  Case of linear topological spaces. Open questions.}

 \vspace{3mm}

 Perhaps, the main result of this report  remains true for the linear topological spaces $  X $ if it is countable normed or metrizable.
In contradiction, in the preprint \cite{Ostrovsky3} is constructed an example of linear topological space $  L  $ and a Borelian r.v. $ \zeta $
with values in this space such that its distribution does not has a compact  embedded linear support of the full measure. \par

\end{document}